\documentstyle[11pt,twoside]{article}
\include{graphicx}
\title{The evaluation of infinite sums of products of two Bessel functions}
\author{\sc R. B.\ Paris \\
{\em Division of Computing and Mathematics}, \\
{\em Abertay University, Dundee DD1 1HG, UK}
}
\begin{document}
\def\f#1#2{\mbox{${\textstyle \frac{#1}{#2}}$}}
\def\dfrac#1#2{\displaystyle{\frac{#1}{#2}}}
\def\boldal{\mbox{\boldmath $\alpha$}}
{\newcommand{\Sgoth}{S\;\!\!\!\!\!/}
\newcommand{\bee}{\begin{equation}}
\newcommand{\ee}{\end{equation}}
\newcommand{\lam}{\lambda}
\newcommand{\ka}{\kappa}
\newcommand{\al}{\alpha}
\newcommand{\fr}{\frac{1}{2}}
\newcommand{\fs}{\f{1}{2}}
\newcommand{\g}{\Gamma}
\newcommand{\br}{\biggr}
\newcommand{\bl}{\biggl}
\newcommand{\ra}{\rightarrow}
\newcommand{\mbint}{\frac{1}{2\pi i}\int_{c-\infty i}^{c+\infty i}}
\newcommand{\mbcint}{\frac{1}{2\pi i}\int_C}
\newcommand{\mboint}{\frac{1}{2\pi i}\int_{-\infty i}^{\infty i}}
\newcommand{\gtwid}{\raisebox{-.8ex}{\mbox{$\stackrel{\textstyle >}{\sim}$}}}
\newcommand{\ltwid}{\raisebox{-.8ex}{\mbox{$\stackrel{\textstyle <}{\sim}$}}}
\renewcommand{\topfraction}{0.9}
\renewcommand{\bottomfraction}{0.9}
\renewcommand{\textfraction}{0.05}
\newcommand{\mcol}{\multicolumn}
\date{}
\maketitle
\pagestyle{myheadings}
\markboth{\hfill \sc R. B.\ Paris  \hfill}
{\hfill \sc  Bessel function sums\hfill}
\begin{abstract}
We examine convergent representations for the sum of Bessel functions
\[\sum_{n=1}^\infty \frac{J_\mu(na) J_\nu(nb)}{n^{\al}}\]
for $\mu$, $\nu\geq0$ and positive values of $a$ and $b$. Such representations enable easy computation of the series in the limit $a, b\to0+$. Particular attention is given to logarithmic cases that occur both when $a=b$ and $a\neq b$ for certain values of $\al$, $\mu$ and $\nu$. 
The series when the first Bessel function is replaced by the modified Bessel function $K_\mu(na)$ is also investigated,
as well as the series with two modified Bessel functions.
\vspace{0.4cm}

\vspace{0.3cm}

\end{abstract}

\vspace{0.3cm}

\noindent $\,$\hrulefill $\,$

\vspace{0.2cm}

\begin{center}
{\bf 1. \  Introduction}
\end{center}
\setcounter{section}{1}
\setcounter{equation}{0}
\renewcommand{\theequation}{\arabic{section}.\arabic{equation}}
We consider the sum\footnote{To avoid overburdening the notation we omit the parameter $\al$ in $S_{\mu,\nu}(a,b)$.}
\begin{equation}\label{e11}
S_{\mu,\nu}(a,b)=\Lambda\sum_{n=1}^\infty \frac{J_\mu(na) J_\nu(nb)}{n^{\al}},\qquad\Lambda=\frac{2^{\mu+\nu}}{a^\mu b^\nu},
\end{equation}
where $J_\nu(z)$ is the Bessel function of the first kind and the multiplicative factor $\Lambda$ is added for convenience.
It is supposed that the orders $\mu$, $\nu\geq 0$ and that $\al$ is a real parameter. It will be further assumed that $a$ and $b$ are positive real quantities. The sum converges absolutely for $\al>0$,
although when $a\neq b$ convergence (non-absolute) is assured when $\al>-1$. Extension of the results to complex $\mu$, $\nu$ and $\alpha$ is straightforward.  

Sums involving the product of $m$ $J$-Bessel functions have been termed $m$-dimensional Schl\"omilch-type series by Miller \cite{ARM}. 
Series of the type in (\ref{e11}) have been encountered in connection with the study of the electromagnetic behaviour of cylindrical antennas in rectangular waveguides. 
When $a=b$, the above sum has been investigated by Williamson \cite{W} under the restriction $\mu+\nu>\al$. This author employed Poisson's summation formula to derive a form suitable for computation when $a$ is small. In the case $\nu=\fs$ and $a=\pi$, it was possible to deduce the value of an infinite sum involving the ${}_4F_3$ hypergeometric series.  More recently, Dominici {\it et al.\/} \cite{DGL}
evaluated (\ref{e11}) in the case $\al-\mu-\nu=-2N$, where $N$ is a non-negative integer. These authors used the representation of the $J$-Bessel function in terms of an integral involving the Gegenbauer polynomial. The case of a sum involving a single $J$-Bessel function has been considered in \cite{TVS}.

The numerical evaluation of (\ref{e11}) becomes difficult in the limit $a, b\to 0+$ on account of the resulting slow convergence of the series. In this paper we obtain a representation of the above series by means of the Mellin transform approach
subject to no additional restrictions on the parameters $\al$, $\mu$ and $\nu$, other than that of the above-mentioned condition on $\al$ for convergence of the series. This approach as a means of dealing with slowly convergent sums was advocated by Macfarlane \cite{M} and is discussed, for example, in the book \cite[Section 4.1.1]{PK}.
We also consider the special case when $\al-\mu-\nu$ is a positive odd integer, both when $a=b$ and $a\neq b$, when logarithmic terms can appear. We conclude with an investigation of the series when one of the $J$-Bessel functions is replaced by a modified $K$-Bessel function, and also series involving two modified Bessel functions.

\vspace{0.6cm}

\begin{center}
{\bf 2. \ The series $S_{\mu,\nu}(a,b)$ when $a=b$}
\end{center}
\setcounter{section}{2}
\setcounter{equation}{0}
\renewcommand{\theequation}{\arabic{section}.\arabic{equation}}
\newtheorem{theorem}{Theorem}
We first discuss the sum in (\ref{e11}) when $a=b$, where we write $S_{\mu,\nu}(a)\equiv S_{\mu,\nu}(a,a)$. Let us introduce the quantity
\[\vartheta:=\al-\mu-\nu;\]
positive odd integer values of $\vartheta$ will be seen to produce logarithmic terms in the expansion of $S_{\mu,\nu}(a)$ and $S_{\mu,\nu}(a,b)$.
Then we have
\begin{equation}\label{e21}
S_{\mu,\nu}(a)=2^\al (\fs a)^\vartheta\sum_{n=1}^\infty \frac{J_\mu(na) J_\nu(na)}{(na)^{\al}}=2^\al  (\fs a)^\vartheta\sum_{n=1}^\infty f(na)\qquad (\al>0),
\end{equation}
where
\[f(x)=\frac{J_\mu(x) J_\nu(x)}{x^{\al}}.\]

From the elementary properties of the Bessel function it is seen that $f(x)=O(x^{-\vartheta})$ as $x\to 0+$ and $f(x)=O(x^{-\al-1})$ as $x\to+\infty$. If we introduce the Mellin transform of $f(x)$ by
\[F(s)=\int_0^\infty t^{s-1} f(t)\,dt=\int_0^\infty \frac{J_\mu(t) J_\nu(t)}{t^\lambda}dt, \qquad \lambda:=\al+1-s\]
valid in the strip of analyticity $\vartheta<\Re (s)<\al+1$, we have by the Mellin inversion theorem (see, for example, \cite[p.~80]{PK})
\begin{equation}\label{e22a}
f(x)=\frac{1}{2\pi i}\int_{c-\infty i}^{c+\infty i} F(s)\,x^{-s}ds,\qquad \vartheta<c<\al+1.
\end{equation}
We observe that, from the convergence condition in (\ref{e21}), the right-hand boundary of the strip of analyticity $\Re (s)=\al+1>1$ and the left-hand boundary is $\Re (s)\leq \al$.
Then we obtain \cite[p.~118]{PK}
\begin{equation}\label{e22}
S_{\mu,\nu}(a)=\frac{2^\al(\fs a)^\vartheta}{2\pi i}\int_{c-\infty i}^{c+\infty i} F(s)\,\zeta(s)\,a^{-s}ds,\qquad \max\{1,\vartheta\}<c<\al+1,
\end{equation}
where $\zeta(s)$ denotes the Riemann zeta function.

From \cite[p.~403]{WBF} we have
\begin{equation}\label{e23}
F(s)=\frac{\g(\lambda)\g(\fs\mu+\fs\nu+\fs-\fs\lambda)}{2^\lambda\g(\fs\lambda+\fs\mu-\fs\nu+\fs)\g(\fs\lambda+\fs\mu+\fs\nu+\fs)\g(\fs\lambda+\fs\nu-\fs\mu+\fs)}
\end{equation}
for $\mu+\nu+1>\Re (\lambda)>0$. It is seen that these conditions correspond to the strip of analyticity in  (\ref{e22a}). The integrand in (\ref{e22}) has simple poles at $s=1$ resulting from $\zeta(s)$ and at
\begin{equation}\label{e23a}
s_m=\vartheta-2m,\qquad m=0, 1, 2, \ldots\ 
\end{equation}
from the numerator Gamma function,
except if $\vartheta$ is a positive odd integer when the pole at $s=1$ is double. There is also a sequence of simple poles on the right of the integration path resulting from $\g(\lambda)$ at $s=\al+1+m$, $m=0, 1, 2, \ldots\ $.

We consider the integral taken round the rectangular contour with vertices at $c\pm iT$, $-d\pm iT$, where $d=2M-\vartheta-1>0$ so that the side parallel to the imaginary axis passes midway between the poles at $s=\vartheta-2M+2$ and $s=\vartheta-2M$.
The contribution from the upper and lower sides $s=\sigma\pm iT$, $-d\leq \sigma\leq c$ as $T\to\infty$ can be estimated by use of the standard results
\begin{equation}\label{e24}
|\g(\sigma\pm it)|\sim \sqrt{2\pi}t^{\sigma-\frac{1}{2}}e^{-\frac{1}{2}\pi t}\qquad (t\to+\infty),
\end{equation}
which follows from Stirling's formula for the gamma function, and \cite[p.~25]{Iv}
\[|\zeta(\sigma\pm it)|=O(t^{{\hat\mu}(\sigma)}\log^\beta t)\qquad (t\to+\infty)\]
where ${\hat\mu}(\sigma)=0$ ($\sigma>1$), $\fs-\fs\sigma$ ($0\leq\sigma\leq1$), $\fs-\sigma$ ($\sigma\leq0$) and $\beta=0$ ($\sigma>1$), 1 ($\sigma\leq 1)$. Then it follows that
\begin{equation}\label{e25}
|F(\sigma\pm it)|=O((\fs t)^{-\frac{1}{2}})\,\frac{\g(\fs(\mu\!+\!\nu\!+\!1)\!-\!\fs(\al\!+\!1\!-\!\sigma)\!+\!\fs it)}{\g(\fs(\mu\!+\!\nu\!+\!1)\!+\!\fs(\al\!+\!1\!-\!\sigma)\!-\!\fs it)}=O((\fs t)^{\sigma-\al-\frac{3}{2}})
\end{equation}
as $t\to +\infty$.
Hence the modulus of the integrand on these horizontal paths is $O(T^{\xi} \log\,T)$ as $T\to\infty$, where $\xi=\sigma+{\hat\mu}(\sigma)-\al-\frac{3}{2}$. Taking into account the different forms of ${\hat\mu}(\sigma)$ and the fact that $\al>0$, we obtain the order estimate $O(T^{-\frac{1}{2}}\log\,T)$ so that the contribution from these paths vanishes as $T\to\infty$.

Displacement of the integration path over the pole at $s=1$, where $\zeta(s)$ has residue 1, and the first $M$ poles of the sequence $\{s_m\}$ we find (provided $\vartheta$ is not a positive odd integer)
\begin{equation}\label{e25b}
S_{\mu,\nu}(a)=2^{\al-1}(\fs a)^{\vartheta-1} F(1)+\sum_{m=0}^{M-1}A_m (\fs a)^{2m}+2^{\mu+\nu} R_M,
\ee
where
\begin{equation}\label{e25a}
A_m=\frac{(-)^m}{m!}\,\frac{\g(1+\mu+\nu+2m)\,\zeta(\vartheta-2m)}{\g(1+\mu+m) \g(1+\nu+m) \g(1+\mu+\nu+m)}
\end{equation}
and the remainder $R_M$ is given by
\begin{equation}\label{e26}
R_M=\frac{a^\vartheta}{2\pi i}\int_{-d-\infty i}^{-d+\infty i}F(s)\,\zeta(s)\,a^{-s}ds
=\frac{a^{2M-1}}{2\pi} \int_{-\infty}^\infty a^{-it} F(-d+it)\zeta(-d+it)\,dt.
\end{equation}

We use the functional relation \cite[p.~603]{DLMF}
\begin{equation}\label{e27}
\zeta(s)=2^s\pi^{s-1} \zeta(1-s) \g(1-s) \sin \fs\pi s,
\end{equation}
together with the fact that $|\zeta(\sigma\pm it)|\leq \zeta(\sigma)$ when $\sigma>1$. Then
\begin{eqnarray}
|\zeta(\vartheta\!-\!2M\!+\!1\!+\!it)|&\leq&\frac{(2\pi)^{\vartheta-2M+1}}{\pi} |\zeta(2M\!-\!\vartheta\!-\!it)|\,|\g(2M\!-\!\vartheta\!-\!it)| \cosh (\fs\pi|t|)\nonumber\\
&\leq& \pi^{\vartheta-2M}\,\zeta(2M\!-\!\vartheta)\,g(t),\label{e27a}
\end{eqnarray}
where
\[g(t)=\pi^{-\frac{1}{2}}|\g(M\!-\!\fs\vartheta\!-\!\fs it) \g(M\!+\!\fs\!-\!\fs\vartheta\!-\!\fs it)|\,\cosh (\fs\pi |t|)\]
\[=O((\fs|t|)^{2M-\vartheta-\frac{1}{2}})\qquad (t\to\pm\infty)\]
upon application of the duplication formula for the gamma function and use of (\ref{e24}).
From (\ref{e25}) we therefore have
\[|F(\vartheta\!-\!2M\!+1\!+\!it)|=O((\fs |t|)^{\vartheta-2M-\al-\frac{1}{2}})\qquad (t\to\pm\infty).\]
Then, since $\zeta(2M-\vartheta)=O(1)$ for large $M$, the modulus of the integrand in (\ref{e26}) is $O((|t|/2)^{-\al-1})$
as $t\to\pm\infty$ and the integral converges ($\al>0$) and is independent of $a$. Consequently we find that  $|R_M|=O((a/\pi)^{2M})$, and hence $R_M\to 0$ as $M\to\infty$ provided $0<a<\pi$.

We therefore see that the upper limit of the summation index on the right-hand side of (\ref{e25b}) can be replaced by $\infty$ provided $0<a<\pi$.
\vspace{0.3cm}

\noindent{\bf 2.1.\ Alternative form of the coefficients $A_m$}
\vspace{0.2cm}

\noindent
If we make use of (\ref{e27}) to express the coefficients $A_m$ in (\ref{e25a}) in terms of a zeta function of positive argument for large $m$, together with the duplication formula for the gamma function, we find
\[A_m=2^{\al-1} \frac{\sin \fs\pi\vartheta}{\pi} \bl(\frac{\pi}{2}\br)^{\vartheta-2m-1} A_m',\]
where
\[A_m'=\frac{\g(\fs\!+\!\fs\mu\!+\!\fs\nu\!+\!m)\g(1\!+\!\fs\mu\!+\!\fs\nu\!+\!m)\g(m\!+\!\fs\!-\!\fs\vartheta)\g(m\!+\!1\!-\!\fs\vartheta)\,\zeta(2m\!+\!1\!-\!\vartheta)}{m! \g(1\!+\!\mu\!+\!m)\g(1\!+\!\nu\!+\!m) \g(1\!+\!\mu\!+\!\nu\!+\!m)}.\]
This yields the expansion in the alternative form (provided $\vartheta$ is not a non-negative integer)
\begin{equation}\label{e29}
S_{\mu,\nu}(a)={\hat F}(1)+2^{\al-1} (\fs\pi)^{\vartheta-1}\frac{\sin \fs\pi\vartheta}{\pi}\sum_{m=0}^\infty A_m' \bl(\frac{a}{\pi}\br)^{2m},
\end{equation}
which was obtained in an equivalent form\footnote{In \cite{W} the factor $\g(\al)$ was omitted in the expression for ${\hat F}(1)$.} in \cite{W}.
In this last reference the condition $\mu+\nu>\al$ (that is, $\vartheta<0$) was imposed; this results in the infinite sequence of poles $s_m$ in (\ref{e23a}) lying entirely in $\Re (s)<0$ so that the poles are all simple. 

Application of the well-known results $\zeta(2m+1-\delta)=O(1)$ and $\g(a+m)/\g(b+m)\sim m^{a-b}$ for large $m$  shows that $A_m'\sim m^{-\al-1}$ as $m\to\infty$. We thus have confirmation that the sums in (\ref{e28}) and (\ref{e29}) converge (since $\al>0$) for the wider domain  $0<a\leq\pi$. It is conjectured that a more refined treatment of the remainder integral $R_M$ in (\ref{e26}), which takes into account the oscillatory nature of the integrand, would produce a more precise estimate that included, in addition to the basic order term $(a/\pi)^{2M}$, a negative power of $M$. This would yield $R_M\to 0$ as $M\to\infty$ for $0<a\leq\pi$.

Then we obtain the following result:

\begin{theorem}$\!\!\!.$\ Let $\vartheta=\al-\mu-\nu \neq 0, 1, 2, \ldots\ $ and $\alpha>0$. Then we have the convergent expansion
\begin{equation}\label{e28}
S_{\mu,\nu}(a)={\hat F}(1)+\sum_{m=0}^\infty A_m (\fs a)^{2m},
\end{equation}
valid for $0<a\leq\pi$, where 
\begin{equation}\label{e28a}
{\hat F}(1)=\frac{(\fs a)^{\vartheta-1} \g(\al)\g(\fs\mu\!+\!\fs\nu\!+\!\fs\!-\!\fs\al)}{2\g(\fs\al\!+\!\fs\mu\!-\!\fs\nu\!+\!\fs)\g(\fs\al\!+\!\fs\mu\!+\!\fs\nu\!+\!\fs)\g(\fs\al\!+\!\fs\nu\!-\!\fs\mu\!+\!\fs)}
\end{equation}
and the coefficients $A_m$ are defined in (\ref{e25a}).
\end{theorem}
\vspace{0.3cm}

\noindent{\bf 2.2.\ The case when $\vartheta$ is a non-negative integer}
\vspace{0.2cm}

\noindent We now consider the case when $\vartheta=N$, where $N=0, 1, 2, \dots\ $. If $\vartheta=2N$, the infinite sequence of poles is $s_m=2N-2m$. Then $\zeta(s_m)$ appearing in the coefficients $A_m$ in (\ref{e25a}) vanishes for $m\geq N+1$ on account of the trivial zeros of $\zeta(s)$ at $s=-2, -4, \ldots\ $. In this case the expression in (\ref{e28}) is modified by the infinite sum being replaced by the sum with index $0\leq m\leq N$; see (\ref{e211a}) below.

When $\vartheta=2N+1$, we have $s_m=2N+1-2m$. There is then a double pole when $m=N$, since the pole at $s_N$ coincides with the pole of $\zeta(s)$ at $s=1$. Letting $s=1+\epsilon$, where $\epsilon\to 0$, we find that the integrand in (\ref{e22}) (including the multiplicative factor $2^\al (\fs a)^\vartheta$) is
\[\frac{(\fs a)^{2N-\epsilon}\zeta(1+\epsilon)\g(\al-\epsilon)\g(-N+\fs\epsilon)}{2\g(N+1+\mu-\fs\epsilon)\g(N+1+\nu-\fs\epsilon)\g(N+1+\mu+\nu-\fs\epsilon)},\]
where $\g(-N+\fs\epsilon)\sim 2\epsilon^{-1}(-)^N/\g(N+1-\fs\epsilon)$. Making use of the results
\[\zeta(1+\epsilon)=\epsilon^{-1}\{1+\epsilon\gamma+O(\epsilon^2)\},\qquad \g(z+\epsilon)=\g(z)\{1+\epsilon \psi(z)+O(\epsilon^2)\},\]
where $\gamma=0.55721\ldots$ is the Euler-Mascheroni constant and $\psi(z)$ is the psi-function, we obtain the expansion of the above integrand about the point $s=1$ given  by
\[\frac{(-)^N(\fs a)^{2N}\g(\al)}{ \g(N+1+\mu)\g(N+1+\nu)\g(N+1+\mu+\nu)N!}\,\frac{1}{\epsilon^2}\bl\{1+e\Upsilon_N(a)+O(\epsilon^2)\br\},\]
where
\begin{equation}\label{e212}
\Upsilon_N(a)=\gamma-\log\,(\fs a)-\psi(\al)+\fs\psi(N+1)+\fs\psi(N+1+\mu)\]\[+\fs\psi(N+1+\nu)+\fs\psi(N+1+\mu+\nu).
\end{equation}
The residue at the double pole is therefore
\[\frac{(-)^N(\fs a)^{2N}\g(\al)\,\Upsilon_N(a)}{\g(N+1+\mu)\g(N+1+\nu)\g(N+1+\mu+\nu)N!}.\]

Hence we have the expansions:
\begin{theorem}$\!\!\!.$\ Let $N$ be a non-negative integer and $\vartheta=\al-\mu-\nu$, with $\al>0$ Then we have the expansions
\begin{equation}\label{e211a}
S_{\mu,\nu}(a)={\hat F}(1)+\sum_{m=0}^{N} A_m (\fs a)^{2m}\qquad (\vartheta=2N)
\end{equation}
and
\begin{equation}\label{e211b}
S_{\mu,\nu}(a)=\mathop{\sum_{m=0}^\infty}_{\scriptstyle m\neq N}A_m (\fs a)^{2m}+\frac{(-)^N(\fs a)^{2N}\g(\al)\,\Upsilon_N(a)}{ \g(N\!+\!1\!+\!\mu)\g(N\!+\!1\!+\!\nu)\g(N\!+\!1\!+\!\mu\!+\!\nu)N!}\qquad (\vartheta=2N+1)
\end{equation}
for $0<a\leq\pi$. The coefficients $A_m$ are defined in (\ref{e25a}), and ${\hat F}(1)$ and $\Upsilon_N(a)$ are given in (\ref{e28a}) and (\ref{e212}).
\end{theorem}

\vspace{0.6cm}

\begin{center}
{\bf 3. \ The series $S_{\mu,\nu}(a,b)$ when $a\neq b$}
\end{center}
\setcounter{section}{3}
\setcounter{equation}{0}
\renewcommand{\theequation}{\arabic{section}.\arabic{equation}}
We now consider the case $a\neq b$ where, without loss of generality, we suppose $a>b$. Following the same procedure described in Section 2, we have the function $f(x)$ and its Mellin transform $F(s)$ given by
\[f(x)=\frac{J_\mu(ax)J_\nu(bx)}{x^{\al}}, \qquad F(s)=\int_0^\infty\frac{J_\mu(at)J_\nu(bt)}{t^\lambda}\,dt
\qquad (\lambda=\al+1-s).\]
The strip of analyticity of the Mellin transform is $\vartheta<\Re (s)<\al+1$, where we recall that $\vartheta=\al-\mu-\nu$. 
From \cite[p.~401]{WBF}
\begin{equation}\label{e32}
F(s)=\frac{(\fs a)^\lambda b^\nu \g(\fs\mu+\fs\nu+\fs-\fs\lambda)}{a^{\nu+1} \g(1+\nu)\g(\fs\mu-\fs\nu+\fs+\fs\lambda)}\times{}_2F_1\bl(\begin{array}{c}\frac{\mu+\nu+1-\lambda}{2}, \frac{\nu-\mu+1-\lambda}{2}\\1+\nu\end{array}\!\!;\frac{b^2}{a^2}\br)
\end{equation}
for $a>b>0$ and $\mu+\nu+1>\Re (\lambda)>-1$, where ${}_2F_1$ denotes the Gauss hypergeometric function. We remark that when $a=b$ the hypergeometric function can be summed by Gauss' theorem (see (\ref{Ags})) to yield the result in (\ref{e23}) subject to the more restrictive condition $\Re (\lambda)>0$.

Then we have
\begin{equation}\label{e32a}
S_{\mu,\nu}(a,b)=\frac{\Lambda}{2\pi i}\int_{c-\infty i}^{c+\infty i} {\tilde F}(s) \zeta(s)(\fs a)^{-s}\,ds,\qquad (\max \{1,\vartheta\}<c<\al+1)
\end{equation}
where ${\tilde F}(s)=(\fs a)^sF(s)$. The poles of the integrand on the left of the integration path are as before, namely at $s=1$ and $s=s_m$, where $s_m$ is defined in (\ref{e23a}). Provided $\vartheta\neq 1, 3, 5, \ldots$ all these poles are simple.

With $s=\sigma\pm it$, we have from (\ref{e24}), (\ref{e32}) and (\ref{A3}) 
\[|{\tilde F}(\sigma\pm it)|= O((\fs t)^{\sigma+\nu-\al-1}) \,\bl|{}_2F_1\bl(\begin{array}{c}\fs(\sigma\!-\!\vartheta)\!\pm\!\fs it,\fs(\sigma\!-\!\vartheta)\!-\!\mu\!\pm\!\fs it\\1+\nu\end{array}\!\!;\chi\br)\br|\]
\bee\label{e33}
=O((\fs t)^{\sigma-\al-\frac{3}{2}})\,\frac{(1+\sqrt{\chi})^{\al+1-\sigma}}{(\sqrt{\chi})^{\nu+\frac{5}{2}}},\qquad \chi:=\frac{b^2}{a^2}\qquad(t\to\infty).
\ee
The order estimate in $t$ is the same as that in (\ref{e25}) for the case $a=b$. Thus the same arguments apply to justify the displacement of the integration path to the left over the first $M$ poles of the sequence $\{s_m\}$. 

A difficulty presents itself with the remainder integral $R_M$ taken along the rectilinear path $\sigma=\vartheta-2M+1$.
It has not been possible to extract the factor\footnote{In the special case $\mu=\nu$, however, a quadratic transformation of the hypergeometric function exists \cite[(15.8.21)]{DLMF} where the extraction of the factor $(1+\sqrt{\chi})^{2M}$ is possible.} $(1+\sqrt{\chi})^{2M}$ from the ${}_2F_1$ function for $t\in (-\infty, \infty)$, which is seen to be present in the above asymptotic estimate. This would indicate that $R_M=
O(((a+b)/(2\pi))^{2M})$ and hence that $R_M\to 0$ as $M\to\infty$ provided that $1<a+b<2\pi$. 

The residue at the pole $s=s_m$ is
\begin{equation}\label{e34c}
\frac{B_m (\fs a)^{2m}}{\g(1+\nu)},\qquad 
B_m=\frac{(-)^m\zeta(\vartheta-2m)}{m! \g(1+\mu+m)}\,{}_2F_1\bl(\begin{array}{c}-m,-m-\mu\\1+\nu\end{array}\!\!;\chi\br).
\end{equation}
The domain of convergence of the infinite sum of these residues can be determined by examining the large-$m$ behaviour of the coefficients $B_m$. From (\ref{e27}) and the properties of the gamma function we find, provided $\vartheta$ is not an even integer,
\[B_m=O(\pi^{-2m}m^{\nu-\al-\frac{1}{2}} )\ {}_2F_1\bl(\begin{array}{c}-m,-m-\mu\\1+\nu\end{array}\!\!;\chi\br)\]
as $m\to\infty$. From (\ref{B1}), the above hypergeometric function possesses the large-$m$ behaviour $O(m^{-\nu-\frac{1}{2}} (1+\sqrt{\chi})^{2m})$ when $0<\chi<1$.
Hence we find
\[B_m (\fs a)^{2m}=O\bl(m^{-\al-1}\bl(\frac{a+b}{2\pi}\br)^{2m} \br)\qquad (m\to\infty),\]
which shows (since $\al>0$) that the sum of the residues $\sum_{m\geq0}B_m(\fs a)^{2m}$ converges when $0<a+b\leq2\pi$
and $\vartheta$ is not an even integer.

Displacement of the integration path in (\ref{e32a}) to the left over the poles at $s=1$ and $s=s_m$, $m\geq 0$ then yields the following result:
\begin{theorem}$\!\!\!.$\ Let $\vartheta=\al-\mu-\nu$ be non-integer, $a>b>0$ and $\alpha>0$. Then we have the convergent expansion
\begin{equation}\label{e34}
S_{\mu,\nu}(a,b)={\hat F}(1)+\frac{1}{\g(1+\nu)} \sum_{m=0}^\infty B_m (\fs a)^{2m}
\end{equation}
for $0<a+b\leq2\pi$, where
\[{\hat F}(1)=\frac{(\fs a)^{\vartheta-1}\g(\fs\mu+\fs\nu+\fs-\fs\al) }{2\g(1+\nu)\g(\fs\mu-\fs\nu+\fs\al+\fs)}\,{}_2F_1\bl(\begin{array}{c}\frac{\mu+\nu+1-\al}{2}, \frac{\nu-\mu+1-\al}{2}\\1+\nu\end{array}\!\!;\frac{b^2}{a^2}\br)\]
and the coefficients $B_m$ are given by (\ref{e34c}).
\end{theorem}

In \cite{DGL}, the value of the parameter $\vartheta$ was taken as $\vartheta=-2N$, $N=0, 1, 2, \ldots\ $, so that $s_m=-2N-2m$. In this case all the terms in the sum in (\ref{e34}) vanish, except when $k=m=0$. Noting that $\zeta(0)=-\fs$, we obtain from (\ref{e34}) 
\[S_{\mu,\nu}(a,b)=\frac{(\fs a)^{-2N-1}\g(N+\fs) }{2\g(1+\nu) \g(\mu-N+\fs)}\ {}_2F_1\bl(\begin{array}{c}N+\fs, N+\fs-\mu\\1+\nu\end{array}\!\!;\frac{b^2}{a^2}\br)\]
\[\hspace{5cm}-\frac{\delta_{N0}}{2\g(1+\mu)\g(1+\nu)},\qquad (\vartheta=-2N, \ N=0, 1, 2, \ldots\ )\]
valid for $a>b$ and $0<a+b\leq2\pi$,
where $\delta_{N0}$ is the Kronecker delta symbol. This is equivalent to the result given in \cite[Theorem 3.1]{DGL}, although there the domain of validity was given as $0<b<a<\pi$. The result when $\alpha=\mu+\nu$ ($\delta=0$) was also considered by Miller \cite[Eq.~(3.5b)]{ARM} who gave the domain of validity as $0<a+b<2\pi$.
\vspace{0.3cm}

\noindent{\bf 3.1.\ The case when $\vartheta$ is a non-negative integer}
\vspace{0.2cm}

\noindent 
The treatment of the case of non-negative integer values of $\vartheta$ follows a similar procedure to that discussed in Section 2.2. When $\vartheta=2N$, $N=0, 1, 2, \ldots$, then $s_m=2N-2m$ and the sum in (\ref{e34}) terminates with the summation index $m$ satisfying $0\leq m\leq N$.

When $\vartheta=2N+1$, then $s_m=2N+1-2m$ and there is a double pole at $s=1$, where the pole $s_N$ coincides with the pole of $\zeta(s)$. With $s=1+\epsilon$, where $\epsilon\to0$, the integrand in (\ref{e32a}) (including the multiplicative factor $\Lambda$) is
\[\frac{(\fs a)^{2N-\epsilon} \zeta(1+\epsilon) \g(-N+\fs\epsilon)}{2\g(1+\nu) \g(N+1+\mu-\fs\epsilon)}\,{}_2F_1\bl(\begin{array}{c}\!\!-N+\fs\epsilon,-N-\mu+\fs\epsilon\\1+\nu\end{array}\!\!;\frac{b^2}{a^2}\br)\]
\[=\frac{(-)^N (\fs a)^{2N}}{\g(1+\nu)\g(N+1+\mu) N!}\,\frac{1}{\epsilon^2} \bl\{1+\epsilon {\hat \Upsilon}_N(a)+O(\epsilon^2)\br\}\,{}_2F_1\bl(\begin{array}{c}\!\!-N\!+\!\fs\epsilon,-N\!-\!\mu\!+\!\fs\epsilon\\1+\nu\end{array}\!\!;\frac{b^2}{a^2}\br),\]
where
\begin{equation}\label{e35c}
{\hat \Upsilon}_N(a)=\gamma-\log\,\fs a+\fs\psi(N+1+\mu)+\fs\psi(N+1).
\end{equation}

The expansion of the hypergeometric function as $\epsilon\to 0$ is discussed in Appendix C. If we let $\chi:=b^2/a^2$ and define
\begin{equation}\label{e35a}
F_m(\mu,\chi):={}_2F_1\bl(\begin{array}{c}-m, -m-\mu\\1+\nu\end{array}\!\!;\chi\br),\qquad (m=0, 1, 2, \ldots\,),
\end{equation}
then from (\ref{C5}) it is found that
\[{}_2F_1\bl(\begin{array}{c}\!\!-N\!+\!\fs\epsilon,-N\!-\!\mu\!+\!\fs\epsilon\\1+\nu\end{array}\!\!;\chi\br)=F_N(\mu,\chi)-\fs\epsilon\Delta_N(\chi)+O(\epsilon^2).\]
The quantity $\Delta_N(\chi)$ is given by
\[\Delta_N(\chi):=\sum_{r=1}^N \bl(\!\!\begin{array}{c}N\\r\end{array}\!\!\br) \bl(\!\!\begin{array}{c}N+\mu\\r\end{array}\!\!\br)\frac{D_r(N,\mu) \,\chi^{r}}{(1+\nu)_r}\hspace{5cm}\]
\begin{equation}\label{e35b}
\hspace{4cm}+\frac{(\mu)_{N+1} \chi^{N+1}}{(1+\nu)_{N+1} (N+1)}\,{}_3F_2\bl(\begin{array}{c}1,1,1-\mu\\N+\nu+2,N+2\end{array}\!\!;\chi\br)
\end{equation}
where, from (\ref{C2a}), the coefficients $D_r(N,\mu)$ are defined by
\begin{eqnarray}
D_r(N,\mu)&:=&r!\sum_{k=0}^{r-1} \bl(\frac{1}{N-k}+\frac{1}{N+\mu-k}\br)\nonumber\\
&=&r!\{\psi(N\!+\!1)+\Psi(N\!+\!1\!+\!\mu)-\psi(N\!+\!1\!-\!r)-\psi(N\!+\!1\!+\!\mu\!-\!r)\}.\label{e35d}
\end{eqnarray}

The residue at the double pole $s=1$ is therefore given by
\[\frac{(-)^N (\fs a)^{2N}}{\g(1+\nu)\g(N+1+\mu) N!}\bl\{{\hat \Upsilon}_N(a) F_N(\mu,\chi)-\fs \Delta_N(\chi)\br\}.\]
Then we have the following theorem.
\begin{theorem}$\!\!\!.$\ Let $N=0, 1, 2, \ldots\,$, $\chi=b^2/a^2<1$ and $0<a+b\leq2\pi$. When $\vartheta=\al-\mu-\nu$ is a non-negative integer the following expansions hold:
\begin{equation}\label{e37a}
S_{\mu,\nu}(a,b)={\hat F}(1)+\frac{1}{\g(1+\nu)} \sum_{m=0}^{N} B_m (\fs a)^{2m}
\end{equation}
when $\vartheta=2N$,
and 
\[S_{\mu,\nu}(a,b)=\frac{(-)^N (\fs a)^{2N}}{\g(1+\nu)\g(N+1+\mu) N!}\bl\{{\hat \Upsilon}_N(a) F_N(\mu,\chi)-\fs \Delta_N(\chi)\br\}\hspace{3cm}\]
\begin{equation}\label{e37b}
\hspace{6cm}+\frac{1}{\g(1+\nu)} \mathop{\sum_{m=0}^\infty}_{\scriptstyle m\neq N}B_m (\fs a)^{2m}
\end{equation}
when $\vartheta=2N+1$.
The coefficients $B_m$ and the quantities ${\hat \Upsilon}_N(a)$, $F_N(\mu,\chi)$ and $\Delta_N(\chi)$ are defined in (\ref{e34c}), (\ref{e35c}), (\ref{e35a}) and (\ref{e35b}). 
\end{theorem}
\bigskip

\vspace{0.3cm}

\noindent{\bf 3.2.\ Two special cases}
\vspace{0.2cm}

\noindent 
We present two cases of the expansion (\ref{e37b}).
When $\vartheta=1$ ($N=0$), we find upon noting that the finite sum appearing in $\Delta_0(\chi)$ is zero, $F_0(\mu,\chi)=1$ and $\psi(1)=-\gamma$ the expansion
\[
S_{\mu,\nu}(a,b)=\frac{1}{\g(1+\nu)\g(1+\mu)}\bl\{\fs\gamma-\log\,\fs a+\fs\psi(1+\mu)-\frac{\mu \chi}{2(1+\nu)}\,{}_3F_2\bl(\begin{array}{c}1,1,1-\mu\\2+\nu, 2\end{array}\!\!;\chi\br)\br\}\]
\begin{equation}\label{e38}
\hspace{7cm} +\frac{1}{\g(1+\nu)} \sum_{m=1}^\infty B_m (\fs a)^{2m} \qquad (\vartheta=0).
\end{equation}
When $\vartheta=3$ ($N=1$), we find with $F_1(\mu,\chi)=1+(1+\mu)\chi/(1+\nu)$ that
\[S_{\mu,\nu}(a,b)=-\frac{(\fs a)^2}{\g(1+\nu)\g(2+\mu)}\bl\{[\fs(1+\gamma)-\log\,\fs a+\fs\psi(2+\mu)]\bl(1+\frac{(1+\mu)\chi}{(1+\nu)}\br)-\fs\Delta_1(\chi)\br\}\]
\begin{equation}\label{e39}
\hspace{6cm}+\frac{1}{\g(1+\nu)} 
\mathop{\sum_{m=0}^\infty}_{\scriptstyle m\neq 1} B_m (\fs a)^{2m},
\end{equation}
where
\[\Delta_1(\chi)=\frac{(2+\mu)\chi}{(1+\nu)}+\frac{\mu(1+\mu)\chi}{2(1+\nu)(2+\nu)}\,{}_3F_2\bl(\begin{array}{c}1,1,1-\mu\\3+\nu,3\end{array}\!\!;\chi\br).\]

When $a=b$ ($\chi=1$), use of the Gauss summation theorem (see (\ref{Ags})) shows that $A_m=B_m/\g(1+\nu)$ ($m\neq N$). From the summations \cite[p.~452]{PBM}
\[{}_3F_2\bl(\begin{array}{c}1,1,1-\mu\\ \nu+\ell,\ell\end{array}\!\!;1\br)=\left\{\begin{array}{ll}
\dfrac{(1+\nu)}{\mu}\{\psi(1+\mu+\nu)-\psi(1+\nu)\} & (\ell=2)\\
\\
\dfrac{2(2+\nu)(2+\mu+\nu)}{\mu(1+\mu)} \{\psi(3+\mu+\nu)-\psi(2+\nu)\} &\\
\hspace{6cm}-\dfrac{2(2+\nu)}{\mu} & (\ell=3),
\end{array}\right.\]
it can be shown after some routine algebra that the expansions (\ref{e38}) and (\ref{e39}) reduce to the result stated in (\ref{e28}) valid for $\chi=1$.

\vspace{0.6cm}
\begin{center}
{\bf 4. \ Sums involving the modified Bessel functions}
\end{center}
\setcounter{section}{4}
\setcounter{equation}{0}
\renewcommand{\theequation}{\arabic{section}.\arabic{equation}}
In this section we investigate two sums involving the modified Bessel functions $K_\nu(x)$ an $I_\nu(x)$. Thus, we consider the sums
\begin{equation}\label{e41}
S_{\mu,\nu}^{(1)}(a,b)=\sum_{n=1}^\infty \frac{K_\mu(an) J_\nu(bn)}{n^\al},\qquad S_{\mu,\nu}^{(2)}(a,b)=\sum_{n=1}^\infty \frac{K_\mu(an) I_\nu(bn)}{n^\al},
\end{equation}
where again we suppose $\mu, \nu\geq0$. In the first sum we require $a>0$, $b>0$ for convergence with $\alpha$ unrestricted. In the second sum we require either $a>b>0$ with $\alpha$ unrestricted, or $a=b>0$ with $\alpha>0$, since $K_\mu(an) I_\nu(bn)\sim (2n\sqrt{ab})^{-1}$ $\exp\,[-(a-b)n]$ as $n\to\infty$. In both cases we shall continue to consider only real values of $\alpha$.

We proceed in the same manner as in Sections 2 and 3. We have the Mellin transform given by \cite[(10.43.26)]{DLMF}
\[F(s)=\int_0^\infty \frac{K_\mu(at) J_\nu(bt)}{t^\lambda}\,dt,\qquad\lambda=1+\al-s\]
\bee\label{e41a}
=\frac{b^\nu (\fs a)^\lambda}{2a^{1+\nu}}\,\frac{\g(\fs\!-\!\fs\lambda\!+\!\fs\nu\!-\!\fs\mu)\g(\fs\!-\!\fs\lambda\!+\!\fs\nu\!+\!\fs\mu)}{\g(1+\nu)}\times {}_2F_1\bl(\begin{array}{c}\frac{1-\lambda+\nu-\mu}{2}, \frac{1-\lambda+\nu+\mu}{2}\\1+\nu\end{array}\!\!;-\frac{b^2}{a^2}\br)
\ee
provided $\nu-\Re (\lambda)\pm\mu>-1$ and $a>0$, $b>0$. Then we obtain
\begin{equation}\label{e42}
S_{\mu,\nu}^{(1)}(a,b)=\frac{1}{2\pi i}\int_{c-\infty i}^{c+\infty i} {\tilde F}(s) \zeta(s) (\fs a)^{-s}ds\qquad (c>\max \{1,\al-\nu\pm\mu\}),
\end{equation}
where ${\tilde F}(s)=(\fs a)^s F(s)$. The poles of the integrand are situated on the left-hand side of the integration path at $s=1$ and at
\begin{equation}\label{e43}
s_m^\pm =\al-\nu\pm\mu-2m,\qquad m=0, 1, 2 \ldots\ .
\end{equation}

Displacement of the integration path to the left over the poles (we omit the details justifying this process) then yields the following result.
\begin{theorem}$\!\!\!.$\ Let $a>0$, $b>0$ and $\mu, \nu\geq0$, with $\alpha$ real but unrestricted and $\chi=b^2/a^2$. Then, provided the poles are all simple, the following expansion holds for $\chi>0$
\[S_{\mu,\nu}^{(1)}(a,b)=F(1)+\frac{b^\nu}{2^{1+\nu}}\sum_{m=0}^\infty \frac{(-)^m \g(\mu-m)\zeta(s_m^+)}{m! \g(1+\nu)}\, F_m(\mu,-\chi)\bl(\frac{a}{2}\br)^{2m-\mu} 
\]
\begin{equation}\label{e46}
+\frac{b^\nu}{2^{1+\nu}}\sum_{m=0}^\infty \frac{(-)^m \g(-\mu-m)}{m! \g(1+\nu)}\,\zeta(s_m^-)F_m(-\mu,-\chi) \bl(\frac{a}{2}\br)^{2m+\mu} 
\end{equation}
where
\[F(1)=\frac{b^\nu(\fs a)^\al}{2a^{1+\nu}}\,\frac{\g(\fs\!-\!\fs\alpha\!+\!\fs\nu\!-\!\fs\mu)\g(\fs\!-\!\fs\alpha\!+\!\fs\nu\!+\!\fs\mu)}{\g(1+\nu)}\hspace{3cm}\]
\[\hspace{5cm}\times {}_2F_1\bl(\begin{array}{c}\frac{1-\alpha+\nu-\mu}{2}, \frac{1-\alpha+\nu+\mu}{2}\\1+\nu\end{array}\!\!;-\chi\br).\]
The functions $F_m(\pm\mu,-\chi)$ are defined in (\ref{e35a}) and the poles $s_m^\pm$ are specified in (\ref{e43}).
\end{theorem}

To determine the domain of convergence of the expansions in (\ref{e46}) we examine the large-$m$ behaviour of the terms. Upon use of the functional relation for $\zeta(s)$ in (\ref{e27}), this behaviour is essentially controlled by
\[T_m=(-)^m\frac{\Gamma(\pm\mu-m)}{m!} (2\pi)^{-2m} \Gamma(1-s_m^\pm)\sin (\fs\pi s_m^\pm) (\fs a)^{2m} F_m(\pm\mu,-\chi)\]
\bee\label{e44a}
=O(m^{-\al+\nu-\frac{1}{2}}) \bl(\frac{a}{2\pi}\br)^{2m} F_m(\pm\mu,-\chi)\qquad (m\to\infty),
\ee
when $\alpha-\nu\pm\mu$ is not an even integer.
From the asymptotic behaviour derived in (\ref{B2}), we have
$F_m(\pm\mu,-\chi)=O(m^{-\nu-\frac{1}{2}}$ $(1+\chi)^m)$ as $m\to\infty$,
so that
\[T_m=O(m^{-\al-1})\,\frac{(a^2+b^2)^m}{(2\pi)^{2m}}\qquad (m\to\infty).\]
Consequently, the expansion in (\ref{e46}) holds in the domain\footnote{If $\al\leq0$ the domain of convergence is $0<\sqrt{a^2+b^2}<2\pi$.}
\bee\label{e45}
0<\sqrt{a^2+b^2}\leq 2\pi\quad (\al>0).
\ee

The treatment of $S_{\mu,\nu}^{(2)}(a,b)$ is similar, since the Mellin transform
\[F(s)=\int_0^\infty \frac{K_\mu(at) I_\nu(bt)}{t^\lambda}\,dt\qquad (a>b>0)\]
is given by (\ref{e41a}) with the argument of the ${}_2F_1$ function replaced by $+b^2/a^2$. Consequently, we obtain the following expansion.
\begin{theorem}$\!\!\!.$\ Let $a\geq b>0$ and $\mu, \nu\geq0$, with $\alpha$ real but unrestricted (if $a>b$) or $\alpha>0$ (if $a=b$) and $\chi=b^2/a^2$. Then, provided the poles are all simple, the following expansion holds for $0<\chi<1$
\[S_{\mu,\nu}^{(2)}(a,b)=F(1)+\frac{b^\nu}{2^{1+\nu}}\sum_{m=0}^\infty \frac{(-)^m \g(\mu-m)\zeta(s_m^+)}{m! \g(1+\nu)}\, F_m(\mu,\chi)\bl(\frac{a}{2}\br)^{2m-\mu}\]
\begin{equation}\label{e44}
+\frac{b^\nu}{2^{1+\nu}}\sum_{m=0}^\infty \frac{(-)^m \g(-\mu-m)}{m! \g(1+\nu)}\,\zeta(s_m^-)F_m(-\mu,\chi) \bl(\frac{a}{2}\br)^{2m+\mu} 
\end{equation}
where
\[F(1)=\frac{b^\nu(\fs a)^\al}{2a^{1+\nu}}\,\frac{\g(\fs\!-\!\fs\alpha\!+\!\fs\nu\!-\!\fs\mu)\g(\fs\!-\!\fs\alpha\!+\!\fs\nu\!+\!\fs\mu)}{\g(1+\nu)}\hspace{3cm}\]
\[\hspace{5cm}\times {}_2F_1\bl(\begin{array}{c}\frac{1-\alpha+\nu-\mu}{2}, \frac{1-\alpha+\nu+\mu}{2}\\1+\nu\end{array}\!\!;\chi\br).\]
The functions $F_m(\pm\mu,\chi)$ are defined in (\ref{e35a}) and the poles $s_m^\pm$ are specified in (\ref{e43}).
\end{theorem}
Following the estimate $T_m$ in (\ref{e44a}), the large-$m$ behaviour of the terms in the expansion in (\ref{e44}) is controlled by
\[T_m=O(m^{-\al+\nu-\frac{1}{2}}) \bl(\frac{a}{2\pi}\br)^{2m}\,F_m(\pm\mu,\chi).\]
From the asymptotic behaviour in (\ref{B1}) we have
$F_m(\pm\mu,\chi)=O(m^{-\nu-\frac{1}{2}}(1+\sqrt{\chi})^{2m})$ as $m\to\infty$ when $0<\chi<1$,
so that
\[T_m=O(m^{-\al-1})\,\bl(\frac{a+b}{2\pi}\br)^{2m}\qquad (m\to\infty).\]
Consequently the expansions in (\ref{e44}) hold in the domain
\bee\label{e47}
0<a+b\leq 2\pi\quad (\al>0)
\ee
with a similar reduced domain $(0,2\pi)$ when $\al\leq 0$; compare (\ref{e45}).
\vspace{0.6cm}

\begin{center}
{\bf 5. \ An example of $S_{\mu,\nu}^{(1)}(a,b)$ when multiple poles are present}
\end{center}
\setcounter{section}{5}
\setcounter{equation}{0}
The expansions presented in (\ref{e46}) and (\ref{e44}) assume all the poles to be simple. Double poles will arise when either (i) one of the poles in the sequences $\{s_m^\pm\}$ coincides with the pole at $s=1$ or (ii) when $\mu=N$, $N=0, 1, 2, \ldots\ $. In this last case the poles $\{s_m^+\}$ are simple for $0\leq m\leq N-1$, with double poles for $m\geq N$. If $\mu=N$ and a double pole from the sequences $\{s_m^\pm\}$ with $m\geq N$ coincides with $s=1$, then there will be a treble pole.

We do not deal with all the cases that can arise here, although the procedure is the routine, albeit laborious, evaluation of residues of the integrand of the appropriate integral. As an illustrative example, we consider the sum $S_{\mu,\nu}^{(1)}(a,b)$ in the case $\mu=2$ and $\al-\nu=3$. This corresponds to the poles $s_m^+=5-2m$ and $s_m^-=1-2m$; the poles at $s=3, 5$ are simple, that at $s=1$ is a treble pole with those at $s=-1, -3, \ldots$ being double poles.  

The residue of the integrand in (\ref{e42}) at the double poles $s=-2m+1$, $m=1, 2, \ldots$ are evaluated in a similar manner to that described in Section 2.2. With $s=-2m+1+\epsilon$, where $\epsilon\to0$, these are given by the coefficient of $\epsilon^{-1}$ in the expansion of
\[\frac{(\fs a)^2 (\fs b)^\nu}{\g(1+\nu)}\,\bl(\frac{a}{4\pi}\br)^{2m-\epsilon}\!\!\frac{(-)^m \zeta(2m-\epsilon) \g(2m-\epsilon)}{\epsilon^2 \g(m\!+\!1\!-\!\fs\epsilon)\g(m\!+\!3\!-\!\fs\epsilon)}\times {}_2F_1\bl(\begin{array}{c}-m\!+\!\fs\epsilon,-m\!-\!2\!+\!\fs\epsilon\\1+\nu\end{array}\!\!;-\chi\br).\]
From (\ref{C5}), the expansion of the hypergeometric function has the form
\[{}_2F_1\bl(\begin{array}{c}-m\!+\!\fs\epsilon,-m\!-\!2\!+\!\fs\epsilon\\1+\nu\end{array}\!\!;-\chi\br)=
F_m(2,-\chi)-\frac{1}{2}\epsilon \Delta_m(-\chi)+O(\epsilon^2),\]
where $F_m(2,-\chi)$ is given in (\ref{e35a}) and, from (\ref{e35b}) when $\mu=2$, 
\[\Delta_m(-\chi)=\sum_{r=1}^m(-)^r\bl(\!\!\begin{array}{c}m\\r\end{array}\!\!\br)\bl(\!\!\begin{array}{c}m+2\\r\end{array}\!\!\br)\frac{D_r(m,2)}{(1+\nu)_r}\chi^r\hspace{4cm}\]
\[\hspace{3cm}+\frac{(-\chi)^{m+1}(m+2)!}{(1+\nu)_{m+1} (m+1)}\bl\{1+\frac{\chi}{(m+2)(m+\nu+2)}\br\}
\]
with the coefficients $D_r(m,2)$ defined by (\ref{e35d}).
This produces the residues of the double poles given by
\bee\label{e51}
\frac{2(\fs a)^2(\fs b)^\nu}{\g(1+\nu)}\,\frac{(-)^m \zeta(2m) \g(2m)}{m! (m+2)!} \bl(\frac{a}{4\pi}\br)^{2m}\bl\{h_m(a) F_m(2,-\chi)-\fs \Delta_m(-\chi)\br\},
\ee
where
\[h_m(a):=\fs\psi(m+1)+\fs\psi(m+3)-\psi(2m)-\frac{\zeta'(2m)}{\zeta(2m)}-\log\,\bl(\frac{a}{4\pi}\br).\]

The residue at the treble pole at $s=1$ is obtained as the coefficient of $\epsilon^{-2}$ in the expansion of
\[\frac{(\fs a)^2 (\fs b)^\nu}{\g(1+\nu)}\,\frac{(a/2)^\epsilon \zeta(1+\epsilon) \g^2(1+\fs\epsilon)}{\epsilon^2(1-\fs\epsilon)(2-\fs\epsilon)}\,{}_2F_1\bl(\begin{array}{c}\fs\epsilon, -2+\fs\epsilon\\1+\nu\end{array}\!\!;-\chi\br).\]
Upon use of the result $\zeta(1+\epsilon)=\epsilon^{-1}\{1+\epsilon\gamma-\epsilon^2\gamma_1+O(\epsilon^3)\}$, where  $\gamma_1=-0.0728158\ldots$ is the first Stieltjes coefficient, and 
\[
{}_2F_1\bl(\begin{array}{c}\fs\epsilon,-2+\fs\epsilon\\1+\nu\end{array}\!\!;-\chi\br)
=1+\frac{1}{2}\epsilon(2-\frac{1}{2}\epsilon)\frac{\chi}{1+\nu}+\frac{1}{2}\epsilon(2-\frac{1}{2}\epsilon)\frac{\chi^2}{(1+\nu)_2 2!}\hspace{2cm}\]
\[\hspace{4.5cm}-\frac{4\epsilon^2\chi^3}{(1+\nu)_3 3!}\bl\{1-\frac{3\cdot 1 \chi}{(4+\nu)4}+\frac{(3)_2(1)_2 \chi^2}{(4+\nu)_2 (4)_2}- \cdots\br\}+O(\epsilon^3)\]
\[=1+\frac{\epsilon\chi}{1+\nu}\bl(1+\frac{\chi}{2(2+\nu)}\br)-\frac{1}{4}\epsilon^2 G(\chi)+O(\epsilon^3),\]
where
\[G(\chi)=\frac{\chi}{1+\nu}\bl(1+\frac{\chi}{2(2+\nu)}\br)+\frac{2\chi^3}{3(1+\nu)_3}\,{}_3F_2\bl(\begin{array}{c}1,1,3\\4+\nu,4\end{array}\!\!;-\chi\br),\]
the residue at $s=1$ is found to be
\[
{\cal R}_1=\frac{(\fs a)^2(\fs b)^\nu}{4\g(1+\nu)} \bl\{\frac{7}{8}-\gamma^2-2\gamma_1-\frac{3}{2}\log\,\fs a+\log^2 \fs a+\frac{\pi^2}{12}\hspace{4cm}\]
\bee\label{e52}
\hspace{4cm}+\frac{2\chi}{1+\nu}\bl(1+\frac{\chi}{2(2+\nu)}\br)\bl(\frac{3}{4}-\log\,\fs a\br)-\frac{1}{2}G(\chi)\br\}.
\ee

Evaluating the terms corresponding to $m=0, 1$ in the first sum in (\ref{e33}), we finally obtain from (\ref{e51}) and (\ref{e52}) the expansion
\[S_{\mu,\nu}^{(1)}(a,b)=\frac{(\fs a)^{-2}(\fs b)^\nu}{2\g(1+\nu)}\bl\{\zeta(5)-\frac{a^2 \zeta(3)}{4}\bl(1+\frac{\chi}{1+\nu}\br)\br\}+{\cal R}_1\]
\bee\label{e53}
+\frac{2(\fs a)^2(\fs b)^\nu}{\g(1+\nu)}\sum_{m=1}^\infty\frac{(-)^m\zeta(2m)\g(2m)}{m! (m+2)!}\,\bl(\frac{a}{4\pi}\br)^{2m}\bl\{h_m(a) F_m(2,-\chi)-\frac{1}{2}\Delta_m(-\chi)\br\}
\ee
valid when $\mu=2$, $\al-\nu=3$ and subject to the condition (\ref{e45}).
\vspace{0.6cm}

\begin{center}
{\bf 6. \ Concluding remarks}
\end{center}
\setcounter{section}{6}
\setcounter{equation}{0}
All the expansions presented have been verified numerically with the aid of {\it Mathematica}. In particular, the large-$m$ behaviour of the terms in the various expansions was examined to verify the domains of convergence given in Theorem 3, (\ref{e45}) and (\ref{e47}). In the computation of the expansion in (\ref{e53}) the term $\zeta'(2m)/\zeta(2m)$ was computed using the command {\tt Zeta}$\,'${\tt [s]/Zeta[s]}; alternatively, the result \cite[(25.2.6)]{DLMF}
\[\frac{\zeta'(2m)}{\zeta(2m)}=\frac{-1}{\zeta(2m)}\sum_{k=2}^\infty \frac{\log\,k}{k^{2m}}\qquad (m=1, 2, \ldots)\]
may be employed.
 
The expansions of the alternating versions of the sums considered in this paper can be deduced from the results of Sections 2--4. For the first sum we have
\[{\hat S}_{\mu,\nu}(a,b)=\Lambda\sum_{n=1}^\infty (-)^{n-1} \frac{J_\mu(an)J_\nu(bn)}{n^\al}=S_{\mu,\nu}(a,b)-2^{1-\vartheta} S_{\mu,\nu}(2a,2b),\]
where we recall that $\Lambda=2^{\mu+\nu}/(a^\mu b^\nu)$ and $\vartheta=\alpha-\mu-\nu$. 
For the alternating sums involving modified Bessel functions
\[{\hat S}_{\mu,\nu}^{(1)}(a,b)=\sum_{n=1}^\infty (-)^{n-1} \frac{K_\mu(an)J_\nu(bn)}{n^\al},\qquad
{\hat S}_{\mu,\nu}^{(2)}(a,b)=\sum_{n=1}^\infty (-)^{n-1} \frac{K_\mu(an)I_\nu(bn)}{n^\al},\]
we have similarly
\[{\hat S}_{\mu,\nu}^{(k)}(a,b)
=S_{\mu,\nu}^{(k)}(a,b)-2^{1-\al} S_{\mu,\nu}^{(k)}(2a,2b)\qquad (k=1, 2).\]
It can be verified that the contribution from the pole of $\zeta(s)$ at $s=1$ is absent
in the expansion of these alternating sums. This is also evident from the integral representations for the alternating sums which take the form
\[\frac{1}{2\pi i}\int_{c-\infty i}^{c+\infty i} F(s) (1-2^{1-s}) \zeta(s) (\fs a)^{-s} ds,\]
where $F(s)$ is the appropriate Mellin transform. The factor $(1-2^{1-s}) \zeta(s)$ appearing in the above integrand is regular at $s=1$.

\vspace{0.6cm}

\begin{center}
{\bf Appendix A: \ The asymptotic behaviour of a Gauss hypergeometric function}
\end{center}
\setcounter{section}{1}
\setcounter{equation}{0}
\renewcommand{\theequation}{\Alph{section}.\arabic{equation}}
Consider the Gauss hypergeometric function
\[{\cal F}(\chi)={}_2F_1\bl(\begin{array}{c}A+i\lambda,B+i\lambda\\C\end{array}\!\!;\chi\br)\qquad (\lambda\to+\infty),\]
where $0<\chi\leq1$. The finite parameters (see Section 3) are $A=\fs(\sigma-\vartheta)$, $B=A-\mu$, $C=1+\nu$ and\footnote{The parameter $\lambda$ in Appendix A is not to be confused with that appearing in Sections 2--4.} $\lambda=\fs t$, with $\sigma<\al+1$. It is easily verified that $C-A-B>0$ and $C-B>0$. 
When $\chi=1$, we have from the Gauss summation theorem 
\bee\label{Ags}
{}_2F_1\bl(\begin{array}{c}a, b\\c\end{array}\!\!;1\br)=\frac{\g(c)\g(c-a-b)}{\g(c-a)\g(c-b)}\qquad(\Re (c-a-b)>0)
\ee
and (\ref{e24}) the result
\[{\cal F}(1)=\frac{\g(C)\g(C\!-\!A\!-\!B\!-\!2i\lambda)}{\g(C\!-\!A\!-\!i\lambda) \g(C\!-\!B\!-\!i\lambda)}\]
\begin{equation}\label{A1}
\sim\frac{\g(C)}{2\sqrt{\pi}}\,2^{C-A-B}\lambda^{-C+\frac{1}{2}}\exp\,\bl[-2i\lambda \log\,2+\frac{\pi iC}{2}-\frac{\pi i}{4}\bl]\quad (\lambda\to+\infty).
\end{equation}

When $0<\chi<1$, we first employ Euler's transformation \cite[(15.8.1)]{DLMF} to yield
\begin{equation}\label{A2}
{\cal F}(\chi)=(1-\chi)^{C-A-B-2i\lambda} {}_2F_1\bl(\begin{array}{c}C\!-\!A\!-\!i\lambda,C\!-\!\!B\!-i\lambda\\C\end{array}\!\!;\chi\br).
\end{equation}
Then, since $C-B>0$ we have the integral representation \cite[(15.6.2)]{DLMF}
\[{}_2F_1\bl(\begin{array}{c}C\!-\!A\!-\!i\lambda,C\!-\!B\!-\!i\lambda\\C\end{array}\!\!;\chi\br)=\frac{\g(C)\g(1\!-\!B\!-\!i\lambda)}{2\pi i \g(C\!-\!B\!-\!i\lambda)} \int_0^{(1+)} h(\tau) e^{-i\lambda\psi(\tau)} d\tau,\] 
where
\[h(\tau)=\frac{\tau^{B-1}(1-\tau)^{C-B-1}}{(1-\chi\tau)^A},\qquad \psi(\tau)=\log\ \bl(\frac{\tau}{(1-\tau)(1-\chi\tau)}\br).\]
The integration path is a loop that starts at $\tau=0$, encircles the point $\tau=1$ in the positive sense (excluding the point $\tau=1/\chi$) and terminates at $\tau=0$. The $\tau$-plane is cut along $(-\infty,1]$ and from the point $1/\chi$ to infinity in a suitable direction.

Stationary points of the phase function $\psi(\tau)$ occur when $\psi'(\tau)=0$; that is, at the points $\tau=\pm 1/\sqrt{\chi}$. The integration path can be deformed to pass over the point $\tau_s=1/\sqrt{\chi}$ in a direction that is locally perpendicular to the real $\tau$-axis. Applying the stationary phase method, where we note that $\psi''(\tau_s)=2\chi^{3/2}/(1-\sqrt{\chi})^2>0$ and make the substitution $\tau-\tau_s=iu$, we have \cite[p.~97]{O}
\[\frac{1}{2\pi i}\int_0^{(1+)} h(\tau) e^{-i\lambda\psi(\tau)} d\tau\sim \frac{h(\tau_s)}{2\pi}e^{-i\lambda\psi(\tau_s)}\int_{-\infty}^\infty e^{i\lambda u^2\psi(\tau_s)}du\]
\[=\frac{h(\tau_s)}{2\pi}e^{-i\lambda\psi(\tau_s)+\pi i/4}\bl(\frac{2\pi}{\lambda \psi''(\tau_s)}\br)^{\!\!1/2}\]
as $\lambda\to+\infty$. 

Since 
\[\psi(\tau_s)=-2\log (1-\sqrt{\chi}), \qquad h(\tau_s)=\frac{\chi^{-C/2}}{(1-\sqrt{\chi})^{C-A-B+1}},\]
we obtain from (\ref{A2}) after some straightforward algebra the estimate 
\[{\cal F}(x)\sim \frac{\g(C)\g(1\!-\!B\!-\!i\lambda)}{2\pi \g(C\!-\!B\!-\!i\lambda)}\,e^{\pi i/4}\sqrt{\frac{\pi}{\lambda}}\,\frac{(1+\sqrt{\chi})^{C-A-B-2i\lambda}}{\chi^{\frac{1}{2}C+\frac{3}{4}}}\]
\begin{equation}\label{A3}
\sim\frac{\g(C)}{2\sqrt{\pi}} \lambda^{-C+\frac{1}{2}}\frac{(1+\sqrt{\chi})^{C-A-B}}{(\sqrt{\chi})^{C+\frac{3}{2}}}\,  \exp\,\bl[-2i\lambda \log (1+\sqrt{\chi})+\frac{\pi iC}{2}-\frac{\pi i}{4}\br]
\end{equation}
for $0<\chi<1$ and $\lambda\to+\infty$.

We remark that if we let $\chi=1$ in (\ref{A3}) then the large-$\lambda$ estimate for ${\cal F}(\chi)_{\chi\to 1}$ agrees with that in (\ref{A1}).
\vspace{0.6cm}

\begin{center}
{\bf Appendix B: \ The asymptotic behaviour of $F_m(\pm\mu, \chi)$ as $m\to\infty$}
\end{center}
\setcounter{section}{2}
\setcounter{equation}{0}
\renewcommand{\theequation}{\Alph{section}.\arabic{equation}}
We consider the asymptotic behaviour for integer $m\to\infty$ of the hypergeometric functions $F_m(\pm\mu,\chi)$ when  (i) $0<\chi<1$ and (ii) $\chi<0$. Application of the transformation \cite[(15.6.2)]{DLMF} shows that
\[F_m(\pm\mu,\chi)\equiv{}_2F_1\bl(\begin{array}{c}-m,-m\mp\mu\\1+\nu\end{array}\!\!;\chi\br)=(1-\chi)^{m\pm\mu}
{}_2F_1\bl(\begin{array}{c}1\!+\!\nu\!+\!m, -m\!\mp\!\mu\\1+\nu\end{array}\!\!;\frac{\chi}{\chi-1}\br).\]
Then, from the expansion given in \cite[(15.12.5)]{DLMF} in terms of the modified Bessel function we obtain, as $m\to\infty$,
\[{}_2F_1\bl(\begin{array}{c}1\!+\!\nu\!+\!m, -m\!\mp\!\mu\\1+\nu\end{array}\!\!;\frac{\chi}{\chi-1}\br)
\sim \frac{\g(1+\nu)}{2\chi^{(1+\nu)/2}} (1-\chi)^{(2+\nu\pm\mu)/2}\,\sqrt{\zeta \sinh \zeta} \ \rho^{-\nu} I_\nu(\rho\zeta)\]
\[\hspace{0.9cm}\sim \frac{\g(1+\nu)}{2\sqrt{\pi}}\,m^{-\nu-\frac{1}{2}}\,\frac{(1-\chi)^{(1+\nu\mp\mu)/2}}{\chi^{\frac{1}{2}\nu+\frac{1}{4}}}
\bl(\frac{1+\sqrt{\chi}}{1-\sqrt{\chi}}\br)^{\rho}\]
where 
\[\rho=m+\fs(1+\nu\pm\mu),\qquad \zeta=\log \bl(\frac{1+\sqrt{\chi}}{1-\sqrt{\chi}}\br),\qquad \sinh \zeta=\frac{2\sqrt{\chi}}{1-\chi}.\]
Hence it follows that
\bee\label{B1}
F_m(\pm\mu,\chi)\sim\frac{\g(1+\nu)}{2\sqrt{\pi}\,m^{\nu+\frac{1}{2}}}\,\frac{(1+\sqrt{\chi})^{2m+1+\nu\pm\mu}}{\chi^{\frac{1}{2}\nu+\frac{1}{4}}}\qquad (0<\chi<1)
\ee
as $m\to\infty$.

For the hypergeometric function with negative argument we have \cite[(18.5.8)]{DLMF}
\[\bl(\!\!\begin{array}{c} m+\nu\\m\end{array}\!\!\br) {}_2F_1\bl(\begin{array}{c}-m,-m\mp\mu\\1+\nu\end{array}\!¬\!;-\chi\br)=(1+\chi)^m P_m^{(\nu,\pm\mu)} \bl(\frac{1-\chi}{1+\chi}\br),\]
where $P_m^{(\al, \beta)}(x)$ is the Jacobi polynomial. We observe that for $-\chi\in (0,-\infty)$ the argument
$(1-\chi)/(1+\chi)\in (-1,1)$. From \cite[(18.15.6)]{DLMF}, we then have the large-$m$ behaviour in terms of the Bessel function
\[P_m^{(\nu,\pm\mu)}\bl(\frac{1-\chi}{1+\chi}\br)\sim \frac{m^\nu}{2^\frac{1}{2} \rho^\nu}\,\frac{\theta^\frac{1}{2} J_\nu(\rho\theta)}{(\sin \fs\theta)^{\nu+\frac{1}{2}} (\cos \fs\theta)^{\pm\mu+\frac{1}{2}}}, \qquad \cos \theta=\frac{1-\chi}{1+\chi}\]
\[\sim\frac{1}{\sqrt{\pi m}}\,\frac{(1+\chi)^{(1+\nu\pm\mu)/2}}{\chi^{\frac{1}{2}\nu+\frac{1}{4}}} \cos\, [2\rho\phi-\fs\pi\nu-\f{1}{4}\pi],\qquad \phi=\arctan \sqrt{\chi}.\]
Hence we obtain the estimate
\bee\label{B2}
F_m(\pm\mu,-\chi)\sim \frac{\g(1+\nu)}{\sqrt{\pi}\,m^{\nu+\frac{1}{2}}}\,
\frac{(1+\chi)^{m+\frac{1}{2}(1+\nu\pm\mu)}}{\chi^{\frac{1}{2}\nu+\frac{1}{4}}} \cos\, [2\rho\phi-\fs\pi\nu-\f{1}{4}\pi]\qquad (\chi>0)
\ee
as $m\to\infty$.

\vspace{0.6cm}

\begin{center}
{\bf Appendix C: \ The small-$\epsilon$ expansion of ${}_2F_1(-N\!+\!\epsilon,-N\!-\!\mu\!+\!\epsilon;1+\nu;\chi)$}
\end{center}
\setcounter{section}{3}
\setcounter{equation}{0}
\renewcommand{\theequation}{\Alph{section}.\arabic{equation}}
Let $N=0, 1, 2, \ldots$, $\mu, \nu\geq 0$, $0<\chi<1$ and $\epsilon$ be a parameter such that $\epsilon\to0$. Then, using the fact that $(-N+\epsilon)_{N+r}=\epsilon (-N+\epsilon)_N(1+\epsilon)_{r-1}$ for $r\geq 1$, we have
\[{}_2F_1\bl(\begin{array}{c}\!\!-N+\epsilon, -N-\mu+\epsilon\\1+\nu\end{array}\!\!;\chi\br)=1+\sum_{r=0}^N\frac{(-N+\epsilon)_r (-N-\mu+\epsilon)_r }{(1+\nu)_r r!}\chi^r\]
\begin{equation}\label{C1}
+\epsilon (-N+\epsilon)_N \sum_{r=1}^\infty\frac{(-N-\mu+\epsilon)_{N+r} (1+\epsilon)_{r-1}}{(1+\nu)_{N+r} (N+r)!}\,\chi^{N+r}.
\end{equation}

We first consider the finite sum and write
\[{\cal S}_1=1+\sum_{r=0}^N\frac{(-N+\epsilon)_r (-N-\mu+\epsilon)_r }{(1+\nu)_r r!}\chi^r.\]
From the expansion
\[(\beta+\epsilon)_r
=(\beta)_r\bl\{1+\epsilon(\psi(\beta+r)-\psi(\beta))+O(\epsilon^2)\br\}
=(\beta)_r\bl\{1+\epsilon\sum_{k=0}^{r-1}\frac{1}{\beta+k}+O(\epsilon^2)\br\},\]
with $\beta$ put equal to $-N$ and $-N-\mu$ in turn, we obtain
\[{\cal S}_1=1+\sum_{r=1}^N\frac{(-N)_r(-N-\mu)_r\chi^r}{(1+\nu)_r r!}\bl\{1-\epsilon \frac{D_r(N,\mu)}{r!}+O(\epsilon^2)\br\}\]
\begin{equation}\label{C2}
=F_N(\mu,\chi)-\epsilon\sum_{r=1}^N \bl(\!\!\begin{array}{c}N\\r\end{array}\!\!\br) \bl(\!\!\begin{array}{c}N+\mu\\r\end{array}\!\!\br)\frac{D_r(N,\mu)}{(1+\nu)_r}\,\chi^r+O(\epsilon^2),
\end{equation}
where $F_N(\mu,\chi)$ is defined in (\ref{e35a}) and
\[D_r(N,\mu):=r!\sum_{k=0}^{r-1}\bl(\frac{1}{N-k}+\frac{1}{N+\mu-k}\br)\]
\begin{equation}\label{C2a}
=r!\{\psi(N\!+\!1)+\psi(N\!+\!1\!+\!\mu)-\psi(N\!+\!1\!-\!r)-\psi(N\!+\!1\!+\!\mu\!-\!r)\}.
\end{equation}

If we denote the infinite sum in (\ref{C1}) by ${\cal S}_2$ then use of the identity $(\beta)_{N+r+1}=(\beta)_{N+1} (\beta+N+1)_r$ shows that
\begin{eqnarray}
{\cal S}_2&=&\epsilon\,(-N+\epsilon)_N \sum_{r=1}^\infty\frac{(-N-\mu+\epsilon)_{N+r} (1+\epsilon)_{r-1}}{(1+\nu)_{N+r} (N+r)!}\,\chi^{N+r}\nonumber\\
&=&\epsilon\,(-N)_N \sum_{r=1}^\infty \frac{(-N-\mu)_{N+r} (r-1)!}{(1+\nu)_{N+r} (N+r)!}\,\chi^{N+r}+O(\epsilon^2)\nonumber\\
&=&\frac{\epsilon\,(-N)_N (-N-\mu)_{N+1}}{(1+\nu)_{N+1} (N+1)!}\,\chi^{N+1} \sum_{r=0}^\infty \frac{(1-\mu)_r r!}{(c+N+1)_r (N+2)_r}\,\chi^r+O(\epsilon^2)\nonumber\\
&=&-\frac{\epsilon\,(\mu)_{N+1} \chi^{N+1}}{(1+\nu)_{N+1} (N+1)}\,{}_3F_2\bl(\begin{array}{c}1,1,1-\mu\\N+\nu+2,N+2\end{array}\!\!;\chi\br)+O(\epsilon^2), \label{C3}
\end{eqnarray}
where the sum has been expressed as a ${}_3F_2$ hypergeometric function.

From (\ref{C2}) and (\ref{C3}) we finally obtain the desired expansion
\begin{equation}\label{C4}
{}_2F_1\bl(\begin{array}{c}\!\!-N\!+\!\epsilon, -N\!-\!\mu\!+\!\epsilon\\1+\nu\end{array}\!\!;\chi\br)=
F_N(\mu,\chi)-\epsilon \Delta_N(\chi)+O(\epsilon^2),
\end{equation}
where
\[\Delta_N(\chi)=\sum_{r=1}^N \bl(\!\!\begin{array}{c}N\\r\end{array}\!\!\br) \bl(\!\!\begin{array}{c}N+\mu\\r\end{array}\!\!\br)\frac{D_r(N,\mu)}{(1+\nu)_r}\,\chi^r\hspace{5cm}\]
\begin{equation}\label{C5}
\hspace{4cm}+\frac{(\mu)_{N+1} \chi^{N+1}}{(1+\nu)_{N+1} (N+1)}\,{}_3F_2\bl(\begin{array}{c}1,1,1-\mu\\N+\nu+2,N+2\end{array}\!\!;\chi\br).
\end{equation}
We remark that when $\mu=0$ the second expression in $\Delta_N(\chi)$ vanishes and that when $\mu=1, 2, \ldots$ the ${}_3F_2(\chi)$ function terminates.

\end{document}